\documentclass{wccm2014}

\usepackage{graphicx}
\usepackage{amsmath}
\usepackage{amsfonts}
\usepackage{amssymb}

\def\enorm#1{|\!|\!| #1 |\!|\!|}

\newcommand{\bn}{\mathbf n}

\newcommand{\bw}{\mathbf w}

\newcommand{\T}{\mathcal T}

\newcommand{\Div}{\mathop{\rm div}}
\newcommand{\DivG}{{\,\operatorname{div_\Gamma}}}
\newcommand{\nablaG}{\nabla_\Gamma}

\newcommand{\Gs}{\mathcal{S}} 

\newcommand{\Wo}{\overset{\circ}{W}}
\newcommand{\la}{\left\langle}
\newcommand{\ra}{\right\rangle}
\title{A SPACE-TIME FEM FOR PDES ON EVOLVING SURFACES}

\author{J\"ORG GRANDE$^{*}$, MAXIM A. OLSHANSKII$^{\dag}$ AND ARNOLD REUSKEN$^{*}$}

\heading{J. Grande, M.A. Olshanskii and A. Reusken}

\address{$^{*}$ Institut f\"ur Geometrie und Praktische  Mathematik\\
RWTH-Aachen University\\
D-52056 Aachen, Germany\\
e-mail: reusken@igpm.rwth-aachen.de
\and
$^{\dag}$Department of Mathematics\\
University of Houston\\
Houston, Texas 77204-3008\\
e-mail: molshan@math.uh.edu
}

\keywords{evolving surface, diffusion, space-time finite elements, discontinuous Galerkin}

\abstract{The paper studies a finite element method for computing transport and diffusion along
evolving surfaces. The method does not require a parametrization of a surface or an extension
of a PDE from a surface into a bulk outer domain. The surface and its evolution may be given
implicitly, e.g., as the solution of a  level set equation. This approach naturally allows
a surface to undergo topological changes and experience local geometric singularities.
The numerical method uses space-time finite elements and is provably second order accurate.
The paper reviews the method, error estimates and shows results for computing
the diffusion of a surfactant on surfaces of two colliding droplets. }

\begin{document}

\section{INTRODUCTION} Partial differential equations posed on evolving surfaces appear in a number of  applications.
Recently, several numerical approaches for handling such type of problems have been introduced, cf.~\cite{DEreview}. In \cite{Dziuk07,DziukElliot2013a} Dziuk and Elliott developed and analyzed a finite element method for computing transport and diffusion on a surface which is based on a \emph{Lagrangian} tracking of the surface evolution. Methods using an \emph{Eulerian} approach were developed in \cite{AS03,DziukElliot2010,XuZh}, based on an extension  of the surface PDE into a bulk domain that contains the surface. Recently, in \cite{ORXsinum,refJoerg,ORXsinum2} another \emph{Eulerian} method, which  does not use an extension of the PDE into the bulk domain, has been introduced and analyzed. The key idea of this method is to
 use restrictions of (usual) space-time volumetric finite element functions to the space-time manifold. This trace finite element technique has been studied for stationary surfaces in  \cite{OlshReusken08,OlshanskiiReusken08,DemlowOlshanskii12}.

In this paper we summarize the key ideas of this space-time trace-FEM and some main results of the error analysis, in particular a result on second order accuracy of the method in space and time. For details we refer to \cite{ORXsinum,ORXsinum2}. In the numerical experiments in \cite{ORXsinum,refJoerg,ORXsinum2} only relatively simple model problems with smoothly evolving surfaces are considered. As a new contribution in this paper we present results of a numerical experiment for a surfactant transport equation on an evolving manifold with a topological singularity, which resembles a droplet collision.
The method that we study uses volumetric finite element spaces which are continuous piecewise linear in space and discontinuous piecewise linear in time. This allows a natural time-marching procedure, in which the numerical approximation is computed on one time slab after another.  Spatial triangulations  may vary per time slab. The results of the numerical experiment show that the method is extremely robust and that even for the case with a topological singularity (droplet collision) accurate results can be obtained on a fixed Eulerian (space-time) grid with a  large time step.

As a model problem we use the following one.
Consider a surface $\Gamma(t)$ passively advected by a \emph{given} smooth velocity field $\bw=\bw(x,t)$, i.e. the normal velocity of $\Gamma(t)$ is given by $\bw \cdot \bn$, with
$\bn$ the unit normal on $\Gamma(t)$. We assume that for all $t \in [0,T] $,  $\Gamma(t)$ is a  hypersurface that is  closed ($\partial \Gamma =\emptyset$), connected, oriented, and contained in a fixed domain $\Omega \subset \Bbb{R}^d$, $d=2,3$. In the remainder we consider $d=3$, but all results have analogs for the case $d=2$. The convection-diffusion equation on the surface that we consider is given by:
\begin{equation}
\dot{u} + ({\Div}_\Gamma\bw)u -{ \nu_d}\Delta_{\Gamma} u= f\qquad\text{on}~~\Gamma(t), ~~t\in (0,T],
\label{transport}
\end{equation}
 with a prescribed source term $f= f(x,t)$ and homogeneous initial condition $u(x,0)=u_0(x)=0$ for $x \in \Gamma_0:=\Gamma(0)$. Here $\dot{u}= \frac{\partial u}{\partial t} + \bw \cdot \nabla u$ denotes the advective material derivative,
 ${\Div}_\Gamma:=\operatorname{tr}\left( (I-\bn\bn^T)\nabla\right)$ is the surface divergence and $\Delta_\Gamma$ is the
 Laplace-Beltrami operator, $\nu_d>0$ is the constant diffusion coefficient. If we take $f=0$ and an initial condition $u_0 \neq 0$, this surface PDE is obtained from mass conservation of the scalar quantity
$u$ with a diffusive flux on $\Gamma(t)$ (cf. \cite{James04,GReusken2011}). A standard transformation to a homogeneous initial condition, which is convenient for a theoretical analysis, leads to \eqref{transport}.
\section{WELL-POSED SPACE-TIME WEAK FORMULATION}
Several weak formulations of \eqref{transport} are known in the literature, see \cite{Dziuk07,GReusken2011}.  The most appropriate for our purposes is a
integral space-time formulation  proposed in \cite{ORXsinum}. In this section we outline this formulation.
Consider  the space-time manifold
 \[
 \Gs= \bigcup\limits_{t \in (0,T)} \Gamma(t) \times \{t\},\quad  \Gs\subset \Bbb{R}^{4}.
 \]
On $L^2(\Gs)$ we use the scalar product $(v,w)_0=\int_0^T \int_{\Gamma(t)} v w \, ds \, dt$.
Let $\nablaG$ denote the tangential gradient for $\Gamma(t)$ and introduce  the space
\[
H=\{\, v \in L^2(\Gs)~|~ \|\nablaG v\|_{L^2(\Gs)} <\infty \, \}
\]
endowed with the scalar product
\begin{equation}
(u,v)_H=(u,v)_0+ (\nablaG u, \nablaG v)_0. \label{inner}
  \end{equation}
We consider the material derivative $\dot{u}$ of $u \in H$ as a distribution on $\Gs$:
\[
  \la\dot u,\phi\ra= - \int_0^T \int_{\Gamma(t)} u \dot \phi + u \phi \DivG \bw\, ds \, dt  \quad \text{for all}~~ \phi \in C_0^1(\Gs).
\]
In \cite{ORXsinum} it is shown that $C_0^1(\Gs)$ is  dense in $H$. If $\dot{u}$ can be extended to a bounded linear functional on $H$, we write $\dot u \in H'$. Define the space
\[
  W= \{ \, u\in H~|~\dot u \in H' \,\}, \quad \text{with}~~\|u\|_W^2 := \|u\|_H^2 +\|\dot u\|_{H'}^2.
\]
In \cite{ORXsinum}   properties of $H$ and $W$ are derived.
Both spaces are Hilbert spaces and smooth functions are  dense in $H$ and $W$.
Define
\[
\Wo:=\{\, v \in W~|~v(\cdot, 0)=0 \quad \text{on}~\Gamma_0\,\}.
\]
The space $\Wo$ is well-defined,
since functions from $W$ have well-defined traces in $L^2(\Gamma(t))$ for any $t\in[0,T]$. 
We introduce the symmetric bilinear form
\[
  a(u,v)= \nu_d (\nablaG u, \nablaG v)_0 + (\DivG \bw\, u,v)_0, \quad u, v \in H,
\]
which is  continuous on $H\times H$:
\begin{equation*}\label{eq:continuity}
a(u,v)\le(\nu_d+\alpha_{\infty}) \|u\|_H\|v\|_H,\quad\text{with}~\alpha_{\infty}:=\|\DivG \bw\|_{L^\infty(\Gs)}.
\end{equation*}
The weak space-time formulation of \eqref{transport} reads: For given $f \in  L^2(\Gs)$ find $u \in \Wo$ such that
\begin{equation} \label{weakformu}
 \la\dot u,v\ra +a (u,v) = (f,v)_0 \quad \text{for all}~~v \in H.
\end{equation}
In \cite{ORXsinum} the inf-sup property
\begin{equation} \label{infsup}
  \inf_{0\neq u \in \Wo}~\sup_{ 0\neq v \in \overset{\phantom{.}}{H}} \frac{\la\dot u,v\ra + a(u,v)}{\|u\|_W\|v\|_H} \geq c_s>0
\end{equation}
is proved. Using this in combination with the continuity result one can show that the \emph{weak formulation \eqref{weakformu} is well-posed.}

We introduce a similar ``\emph{time-discontinuous}'' weak formulation that is better suited for the finite element method that we consider.
 We take a  partitioning of the time interval:  $0=t_0 <t_1 < \ldots < t_N=T$, with a uniform time step $\Delta t = T/N$. The assumption of a uniform time step is made to simplify the presentation, but is not essential. A time interval is denoted by $I_n:=(t_{n-1},t_n]$.  The symbol $\Gs^n$ denotes the space-time interface corresponding to $I_n$, i.e.,  $\Gs^n:=\cup_{t \in I_n}\Gamma(t)\times\{t\}$, and $\Gs:= \cup_{1 \leq n \leq N}  \Gs^n $. We introduce the following subspaces of $H$:
\[H_n:=\{\, v \in H~|~v=0  \quad \text{on}~~\Gs \setminus \Gs^n\, \},
\]
and define the spaces
\begin{align}
  W_n & = \{ \, v\in H_n~|~\dot v \in H_n' \,\}, \quad \|v\|_{W_n}^2 = \|v\|_{H}^2 +\|\dot v\|_{H_n'}^2, \\
W^b  &  := \oplus_{n=1}^N W_n,~~\text{with norm}~~ \|v\|_{W^b}^2= \sum_{n=1}^N \|v\|_{W_n}^2. \label{brokenW}
\end{align}
For $u \in W_n$, the one-sided  limits
$u_+^{n}=u_+(\cdot,t_{n})$ (i.e., $t \downarrow t_n$) and ${u}_{-}^n=u_{-}(\cdot,t_n)$ (i.e., $t \uparrow t_n$) are well-defined in $L^2(\Gamma(t_n))$.  At $t_0$ and $t_N$  only $u_+^{0}$ and $u_{-}^{N}$ are defined.
For $v \in W^b$, a jump operator  is defined by $[v]^n= v_+^n-v_{-}^n \in L^2(\Gamma(t_n))$, $n=1,\dots,N-1$. For $n=0$, we define $[v]^0=v_+^0$.
On the cross sections $\Gamma(t_n)$, $0 \leq n \leq N$,  of $\Gs$ the $L^2$ scalar product is denoted by
\[
 (\psi,\phi)_{t_n}:= \int_{\Gamma(t_n)} \psi \phi \, ds .
\]
In addition to $a(\cdot,\cdot)$,  we define on the broken space $W^b$ the following bilinear forms:
\begin{align*}
 d(u,v)  =  \sum_{n=1}^N d^n(u,v), \quad d^n(u,v)=([u]^{n-1},v_+^{n-1})_{t_{n-1}},\quad
 \la \dot u ,v\ra_b =\sum_{n=1}^N  \la \dot u_n, v_n\ra.
\end{align*}
One can show  that the unique solution to \eqref{weakformu} is also the unique solution of the  following variational problem in the broken space: Find $u \in W^b$ such that
\begin{equation} \label{brokenweakformu}
  \la \dot u ,v\ra_b +a(u,v)+d(u,v) =( f,v )_0 \quad \text{for all}~~v \in W^b.
\end{equation}
For this time discontinuous weak formulation an inf-sup stability result (that is weaker than the one in \eqref{infsup}) can be derived.
The variational formulation uses $W^b$, instead of $H$,  as test space, since the term $d(u,v)$ is not well-defined for an arbitrary $v\in H$.
Also note that the initial condition $u(\cdot,0)=0$ is not an essential condition in the space $W^b$ but is treated in a weak sense (as is standard in DG methods for time dependent problems).
 From an algorithmic point of view  the formulation \eqref{brokenweakformu} has the advantage that due to the use of the broken space $W^b= \oplus_{n=1}^N W_n$ it can be solved in a time stepping manner.
\section{SPACE-TIME FINITE ELEMENT METHOD } We introduce a finite element method which is a Galerkin method with $W_h \subset  W^b$ applied  to the variational formulation \eqref{brokenweakformu}.
To define this $W_{h}$, consider  the partitioning of the space-time volume domain $Q= \Omega \times (0,T] \subset \Bbb{R}^{3+1}$ into time slabs  $Q_n:= \Omega \times I_n$. Corresponding to each time interval $I_n:=(t_{n-1},t_n]$ we assume a given shape regular tetrahedral triangulation $\T_n$ of the spatial domain $\Omega$. The corresponding spatial mesh size parameter is denoted by $h$. 
Then $\mathcal{Q}_h=\bigcup\limits_{n=1,\dots,N}\T_n\times I_n$ is a subdivision of $Q$ into space-time
prismatic nonintersecting elements. We shall call $\mathcal{Q}_h$ a space-time triangulation of $Q$.
Note that this triangulation is not necessarily fitted to the surface $\Gs$. We allow $\T_n$ to vary with $n$ (in practice, during time integration one may wish to adapt the space triangulation depending on the changing local geometric properties of the surface) and so
the elements of $\mathcal{Q}_h$ may not  match at $t=t_n$.

For any $n\in\{1,\dots,N\}$, let $V_n$ be the finite element space of continuous piecewise linear functions on $\T_n$.
We define the \emph{volume space-time finite element space}:
\begin{equation} \label{defVn}
V_{h}:= \{ \, v: Q \to \Bbb{R} ~|~  v(x,t)= \phi_0(x) + t \phi_1(x)~\text{on every}~Q_n,~\text{with}~\phi_0,\, \phi_1 \in V_n\,\}.
\end{equation}
Thus, $V_{h}$ is a space of piecewise P1 functions with respect to  $\mathcal{Q}_h$, continuous in space and discontinuous in time. Now we define our \emph{surface finite element space} as the space
of traces of functions from $V_{h}$ on $\Gs$:
\begin{equation} \label{deftraceFE}
  W_{h} := \{ \, w:\Gs \to \Bbb{R}~|~ w=v_{|\Gs}, ~~v \in V_{h} \, \}.
\end{equation}
The finite element method reads: Find $u_h \in W_{h}$ such that
\begin{equation} \label{brokenweakformu_h}
  \la \dot u_h ,v_h\ra_b +a(u_h,v_h)+d(u_h,v_h) = (f,v_h)_0 \quad \text{for all}~~v_h \in W_{h}.
\end{equation}
As usual in time-DG methods,  the initial condition for $u_h(\cdot,0)$ is treated in a weak sense.
Due to $u_h\in H^1(Q_n)$ for all $n=1,\dots,N$, the first term  in \eqref{brokenweakformu_h}
can be written as
\[
\la \dot u_h ,v_h\ra_b=\sum_{n=1}^N\int_{t_{n-1}}^{t_n}\int_{\Gamma(t)} (\frac{\partial u_h}{\partial t} +\bw\cdot\nabla u_h)v_h ds\,dt.
\]
The method can be implemented with a time marching strategy. Of course, for the implementation of the method one needs a quadrature rule to approximate the integrals over $\Gs^n$. This issue is briefly addressed in Section~\ref{sectexp}.
\section{DISCRETIZATION ERROR ANALYSIS}
In this section we briefly address  the discretization error analysis of the method \eqref{brokenweakformu_h}, which is presented in  \cite{ORXsinum2}. We first explain a discrete mass conservation property of the scheme \eqref{brokenweakformu_h}. We consider the case that \eqref{transport}  is derived from mass conservation of a scalar quantity with a diffusive flux on $\Gamma(t)$. The original problem then has a nonzero initial condition $u_0$ and a source term $f \equiv 0$.  The solution $u$ of the original problem has the mass conservation property $\bar u(t):=\int_{\Gamma(t)} u \, ds =\int_{\Gamma(0)} u_0 \, ds $ for all $t \in [0,T]$. After a suitable transformation one obtains the equation   \eqref{transport} with a zero initial condition $u_0$ and a right hand-side $f$ which satisfies $\int_{\Gamma(t)} f \, ds =0 $ for all $t \in [0,T]$. The solution $u$ of \eqref{transport} then has the ``shifted'' mass conservation property $\bar u(t)=0$ for all $t \in [0,T]$.  Tak ing suitable test functions in the discrete problem \eqref{brokenweakformu_h} we obtain that the discrete solution $u_h$ has the following weaker mass conservation property, with $\bar u_h(t):=\int_{\Gamma(t)} u_h \, ds$:
\begin{equation}\label{means}
\bar u_{h,-}(t_{n})=0\quad\text{and}\quad\int_{t_{n-1}}^{t_n} \bar u_h(t)\, dt=0,\quad n=1,2,\dots N.
\end{equation}
For a \emph{stationary} surface, $\bar u_h(t)$ is a piecewise affine function and  thus \eqref{means} implies $\bar u_h(t)\equiv0$, i.e,. we have exact mass conservation on the discrete level.
 If the surface evolves, the finite element method is not necessarily mass conserving:  \eqref{means} holds, but $\bar u_h(t) \neq 0$ may occur for $ t_{n-1} \le t < t_n$.  In the discretization error analysis we use a \emph{consistent} stabilizing term involving the quantity $\bar u_h (t)$.  More precisely, define
\begin{equation} \label{stabblf}
  a_\sigma(u,v):= a(u,v)+\sigma \int_0^T \bar u(t) \bar v(t) \, dt, \quad \sigma \geq 0.
\end{equation}
Instead of \eqref{brokenweakformu_h} we consider the stabilized version: Find $u_h \in W_{h}$ such that
\begin{equation} \label{brokenweakformu_h1}
  \la \dot u_h ,v_h\ra_b +a_\sigma(u_h,v_h)+d(u_h,v_h) =(f,v_h)_0 \quad \text{for all}~~v_h \in W_{h}.
\end{equation}
Taking $\sigma >0$ we expect both a stabilizing effect and an improved discrete mass conservation property. Ellipticity
 of finite element method bilinear form and error bounds are derived in the  mesh-dependent norm:
\[
\enorm{u}_h:=\left(\|u_{-}^N\|_T^2 + \sum_{n=1}^N \|[u]^{n-1}\|_{t_{n-1}}^2+\|u\|_H^2\right)^\frac12.
\]
In the error analysis we need a condition which plays a similar role as the  condition ``$ c - \frac12 \Div b >0$''  used in standard analyses of variational formulations of the convection-diffusion equation $- \Delta u + b\cdot \nabla u + c u=f$ in an Euclidean domain $\Omega \subset \Bbb{R}^n$, cf. \cite{TobiskaBook}. This condition is as follows:
there exists a $c_0 >0$ such that
\begin{equation} \label{ass7}
   \DivG \bw (x,t) +\nu_d c_F(t) \geq c_0 \quad \text{for all}~~x \in \Gamma(t),~t \in [0,T].
\end{equation}
Here $ c_F(t) >0$ results from the Poincare inequality
\begin{equation}\label{Fr}
\int_{\Gamma(t)} |\nabla_\Gamma u|^2 \, ds \geq c_F(t) \int_{\Gamma(t)} ( u- \frac{1}{|\Gamma(t)|} \bar u)^2 \, ds\quad \forall~t\in[0,T],~~\forall~u \in H.
\end{equation}

A main result derived in \cite{ORXsinum2} is given in the following theorem. We assume that the time step $\Delta t$ and the spatial mesh size parameter $h$ have comparable size: $\Delta t \sim h$. \\
{\bf Theorem 1.} \emph{Assume \eqref{ass7}  and take $ \sigma \geq  \frac{\nu_d}2\max\limits_{t\in[0,T]}\frac{c_F(t)}{|\Gamma(t)|}$, where $c_F(t)$ is defined in \eqref{Fr}. Then the ellipticity estimate
\begin{equation} \label{coer}
 \la \dot u,u\ra_b +a_\sigma(u,u)+d(u,u) \geq c_s \enorm{u}_h^2 \quad \text{for all}~~u \in W^b
\end{equation}
holds, with $c_s=\frac12\min\{1,\nu_d,c_0\}$ and $c_0$ from \eqref{ass7}. Let $u \in \Wo$ be the solution of \eqref{weakformu} and assume $u \in  H^2(\Gs)$.  For  the  solution  $u_h \in W_h$  of the discrete problem \eqref{brokenweakformu_h1} the following error bound holds:
\[
\enorm{u-u_h}_h \leq c h \|u\|_{H^2(\Gs)}.
\]
}
\\
A further main result  derived in \cite{ORXsinum2} is related to  second order convergence.  Denote by $\|\cdot\|_{-1}$ the norm dual to the $H^1_0(\Gs)$ norm with respect to the $L^2$-duality. Under the conditions given in Theorem~1 and some further mild assumptions the  error bound
\[
 \|u-u_h\|_{-1} \leq c h^2 \|u\|_{H^2(\Gs)}
\]
holds. This second order convergence is derived in a norm weaker than the commonly considered
$L^2(\Gs)$ norm. The reason  is that our arguments use isotropic  polynomial interpolation error
bounds on 4D space-time elements. Naturally, such bounds call for isotropic space-time $H^2$-regularity bounds for the  solution.
For our problem class such regularity is more restrictive than in an elliptic case, since the solution is generally less regular in time than in space. We can  overcome this by measuring the error in the weaker $\|\cdot\|_{-1}$-norm.

\section{NUMERICAL EXPERIMENT} \label{sectexp}
In \cite{ORXsinum,refJoerg} results of numerical experiments are presented. The examples considered there have smoothly evolving surfaces (e.g. a shrinking sphere) and the results show a convergence of order 1 in an $L^2(H^1)$-norm (i.e. $L^2$ w.r.t time and $H^1$ w.r.t space) and of order 2 in an $L^\infty(L^2)$ norm. This convergence behavior occurs already on relatively coarse meshes and there is no (CFL-type) condition on $\Delta t$.

In the example in this paper we consider an evolving surface $\Gamma(t)$ which undergoes a change of topology and experiences a local singularity.
The computational domain is $x \in \Omega=(-3,3)\times(-2,2)^2$, $t \in [0,1]$.   For representation of the evolving surface we use a level set function
$\phi$ defined as:
\[
  \phi(x,t) = 1 - \frac{1}{\| x -c_+(t)\|^3} - \frac{1}{\| x -c_-(t)\|^3},
\]
with  $c_\pm(t)= \pm\frac32(t - 1, 0 , 0)^T$. The surface $\Gamma(t)$ is defined as the zero level of $\phi(x,t)$, $t \in [0,1]$. Take $t=0$. Then
 for $x \in B(c_+(0);1)$ we have $\| x -c_+(0)\|^{-3} =1$ and $\|x -c_-(0)\|^{-3} \ll 1$. For  $ x \in B(c_-(0);1)$ we have $\| x -c_+(0)\|^{-3}  \ll 1$ and $\| x -c_-(0)\|^{-3} =1$. Hence,  the initial configuration $\Gamma(0)$ is (very) close to two balls of radius $1$, centered at $\pm(1.5, 0, 0)^T$. For $t=1$ the surface $\Gamma(1)$ is the ball around $0$ with radius $2^{1/3}$. For $t >0 $ the two spheres approach each other until time $\tilde t= 1-\tfrac23 2^{1/3}\approx 0.160$, when they touch at the origin. For $t \in (\tilde t,1]$ the surface $\Gamma(t)$ is simply connected and smoothly deforms into the sphere $\Gamma(1)$.

\begin{figure}[ht!]
  \centering
  \includegraphics[width=0.7\textwidth]{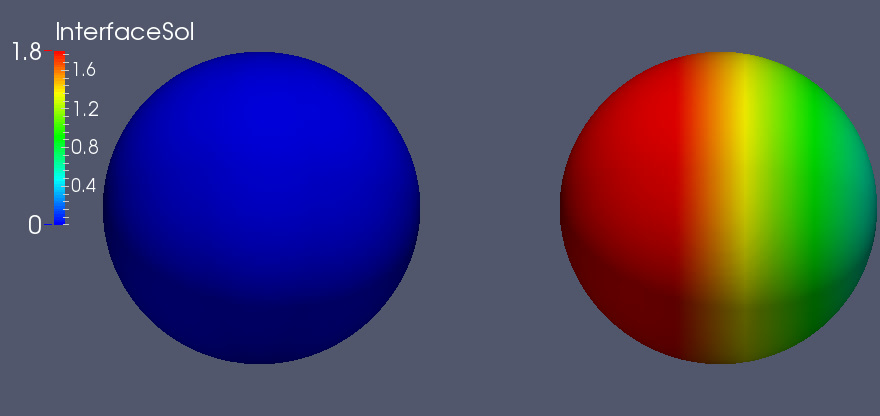}
  \caption{\label{fig:collision-initial}Initial condition as color on the initial zero level $\Gamma(0)$.}
\end{figure}

In the vicinity of $\Gamma(t)$, the gradient $\nabla\phi$ and the time derivative $\partial_t\phi$ are well-defined and given by simple algebraic expressions. We construct the normal wind field, which transports $\Gamma (t)$, by inserting the ansatz $\bw(x,t)=\alpha(x,t)\nabla\phi(x,t)$ into the level set equation $\partial_t\phi + \bw \cdot\nabla\phi=0$. This yields
\[
  \bw= -\frac{\partial_t\phi}{\lvert\nabla\phi\rvert^2}\nabla\phi.
\]
We consider the surfactant advection-diffusion equation
\begin{equation} \label{surfi}
  \begin{cases}
    \dot{u} + \operatorname{div}_\Gamma \bw\, u - \Delta_\Gamma u= 0 & \text{on }\Gamma(t),~t\in(0, 1],\\
    u(\cdot, 0)= u_0 &  \text{on }\Gamma(0).
  \end{cases}
\end{equation}
The initial surfactant distribution is given by
\[
  u_0(x)= \begin{cases} 3 - x_1 & \text{for }x_1 \ge 0,\\
                        0 & \text{else}.
          \end{cases}
\]
The initial configuration is illustrated in Figure~\ref{fig:collision-initial}.

For the construction of a volume space-time finite element space we proceed as follows. On $\Omega$ we start with a level $l=0$
Kuhn-triangulation with mesh width $h_0=2$. We use regular refinement in the vicinity of the interface $\Gamma(t)$ to ensure that the interface is embedded in tetrahedra with refinement level $l \geq 1$. These tetrahedra have the mesh width $h_l=2^{1-l}$.
On each time slab  a level $l$ triangulation is used to define the volume space-time finite element space as in \eqref{defVn}. For simplicity we use the same value for $l$ on all time slabs. The outer space induces a surface finite element space $W_h$ as in \eqref{deftraceFE}. This space is used for a Galerkin discretization of \eqref{surfi}, as given in \eqref{brokenweakformu_h} (note that in this experiment we take $f=0$ and a nonhomogeneous initial condition $u_0$).

We outline the quadrature method for the approximation of integrals over $\Gs^n$. More details are given in \cite{refJoerg}. Consider a single space-time prism $T \times I_n$ that is intersected by  the space time manifold $\Gs^n$. Here $T$ is a level $l$ tetrahedron from the spatial triangulation.  First $T$ is regularly refined into 8 tetrahedra $T_j$, $j=1,\ldots,8$. Each of the resulting space-time prisms $T_j \times I_n$ is partitioned into 4 pentatopes by inserting adequate diagonals. On each of these pentatopes the \emph{linear interpolant}  (in $\Bbb{R}^4$) of the level set function $\phi$ is computed. The zero level of this interpolant is (if not degenerated) a 3-dimensional convex polytope, which can be partitioned into tetrahedra. On these tetrahedra standard quadrature rules can be used. Note that in this approximation procedure we have a \emph{geometric error} due to the approximation of the zero level of $\phi$ (which is the surface) by the zero level of its linear interpolant. This assembling procedure is completely local and can be done prism per prism.

\begin{figure}[ht!]
 \begin{centering}
\begin{minipage}{0.48\textwidth}
  \includegraphics[width=\textwidth]{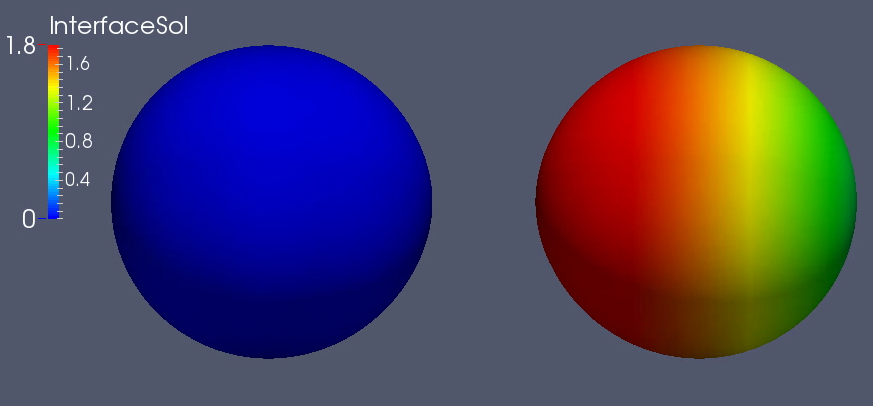}
 \end{minipage} \hfill
\begin{minipage}{0.48\textwidth}
  \includegraphics[width=\textwidth]{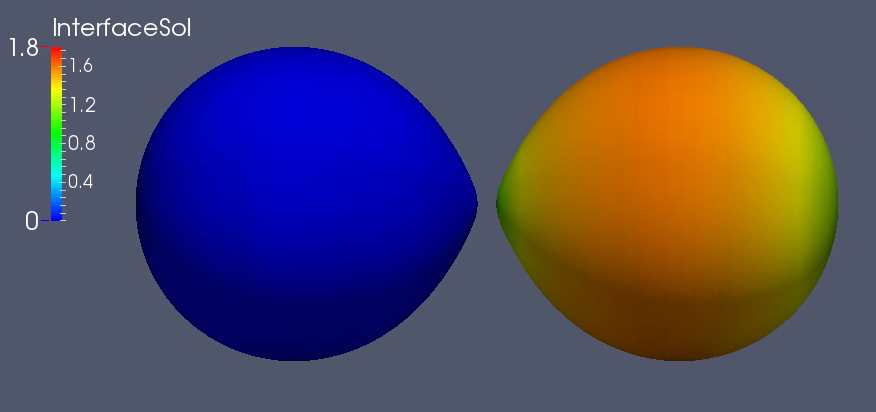}
  \end{minipage}\\[1ex]
\begin{minipage}{0.48\textwidth}
  \includegraphics[width=\textwidth]{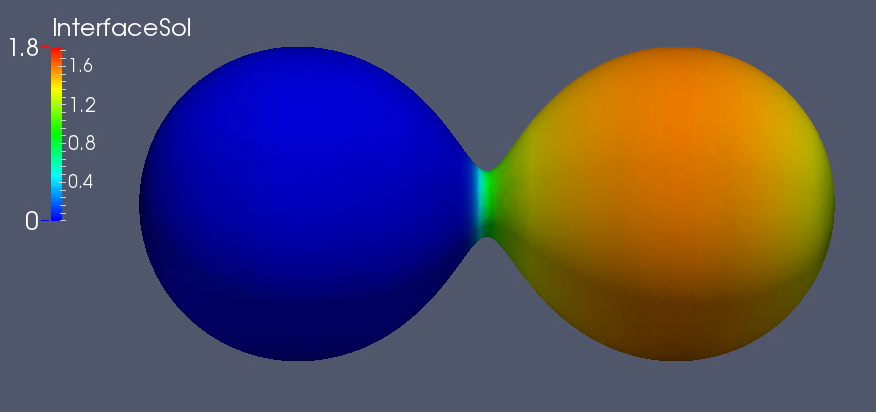}
 \end{minipage} \hfill
\begin{minipage}{0.48\textwidth}
  \includegraphics[width=\textwidth]{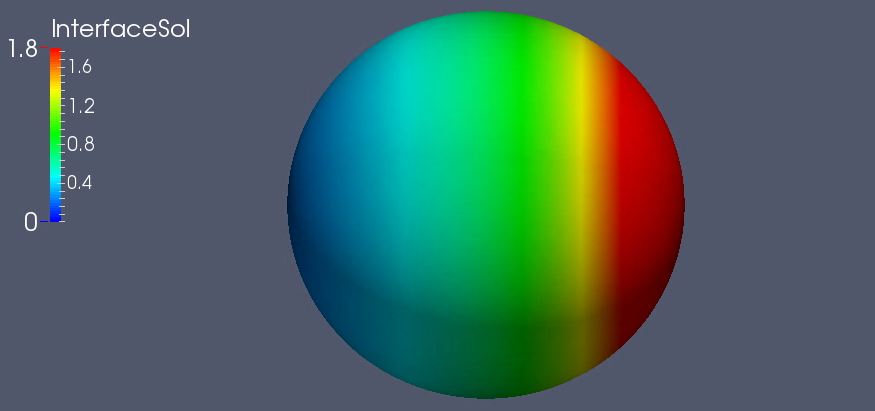}
  \end{minipage}
  \caption{\label{fig3} Snapshots of discrete solution, $l=5$, $\Delta t=2^{-7}$.}
\end{centering}
\end{figure}

\begin{figure}[ht!]
  \centering
  \includegraphics[width=0.7\textwidth]{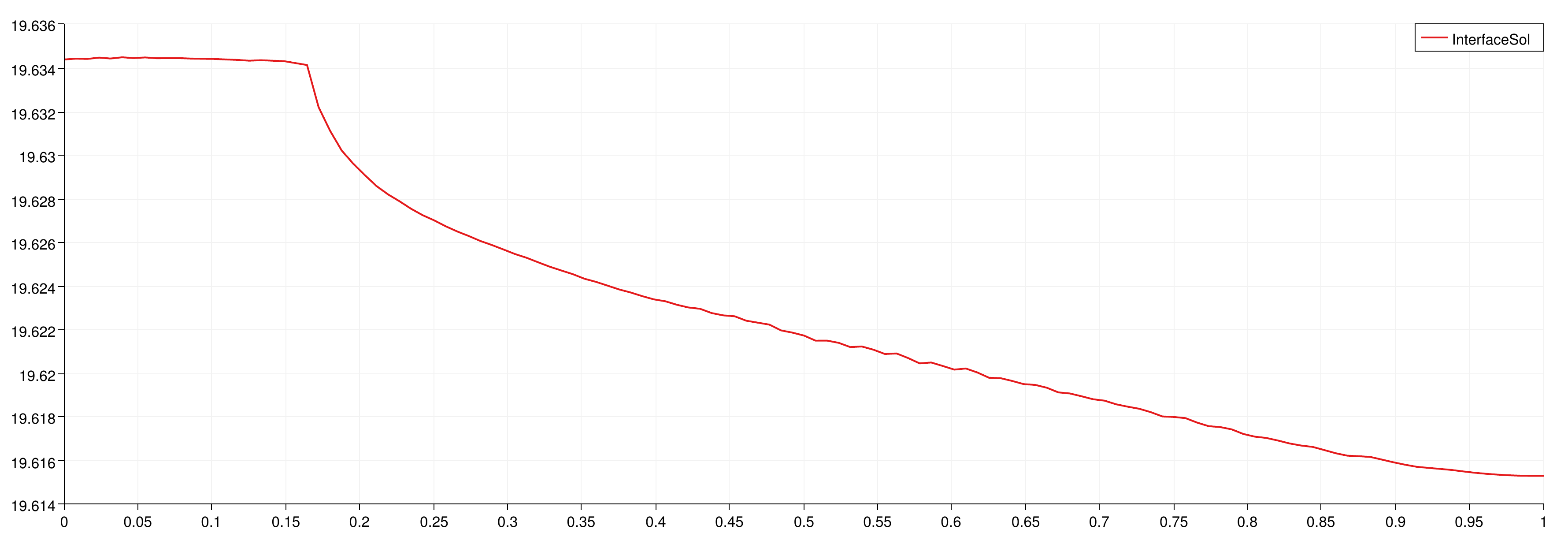}
  \caption{\label{fig:collision-mass-fine}The total amount of surfactant $I_{5,1/128}$ over time.}
\end{figure}

We present some results of numerical experiments. In Figure~\ref{fig3} we show a few snapshots of the surface and the computed surfactant distribution on a relatively fine space-time mesh, namely level $l=5$ and $\Delta t =2^{-7}$. As a measure of accuracy we computed the discrete mass on the space-time manifold:
\[
 I_{l,dt}(t_n)= \int_{\Gamma_h(t_n)} u\,d\sigma,\quad n=0,1,\ldots, N.
\]
where $\Gamma_h(t_n)$ is the approximation of $\Gamma(t_n)$ obtained as zero level of the piecewise linear interpolant of $\phi$, cf. explanation  above.
\\
For $l=5$, $\Delta t=1/128$ the result is shown in  Figure~\ref{fig:collision-mass-fine}. We interpolated the values  $I_{l,dt}(t_n)$, $n=0, \ldots, N$, resulting in the discrete mass quantity as a function of $t \in [0,1]$. There is a mass loss of about $0.018$, which corresponds to a relative error of $\sim~ 9\cdot 10^{-4}$.

\begin{figure}[ht!]
  \centering
  \includegraphics[width=0.7\textwidth]{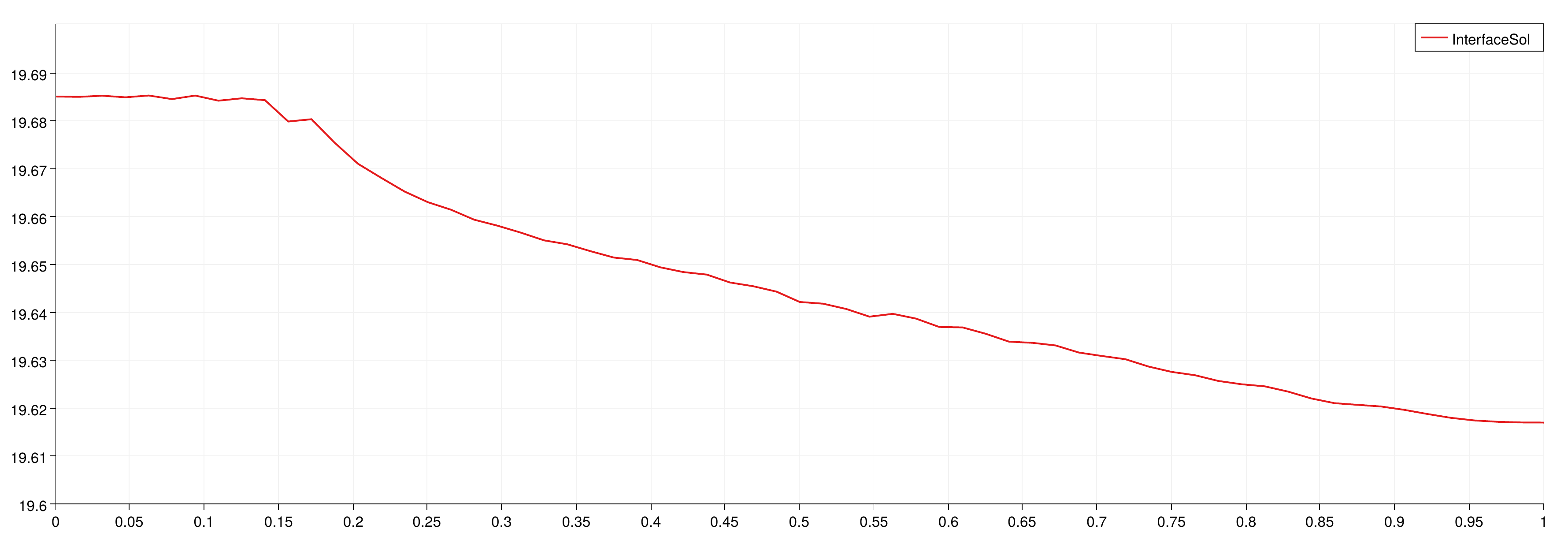}
  \caption{\label{fig:collision-mass}The total amount of surfactant $I_{4,1/64}$ over time.}
\end{figure}

In  Figure~\ref{fig:collision-mass} we show the result for $l=4$, $\Delta t=1/64$. The mass loss is $\sim~0.065$, which is about a factor $3.6$ more than for the case $l=5$, $\Delta t=1/128$.

\begin{figure}[ht!]
 \begin{centering}
\begin{minipage}{0.48\textwidth}
  \includegraphics[width=\textwidth]{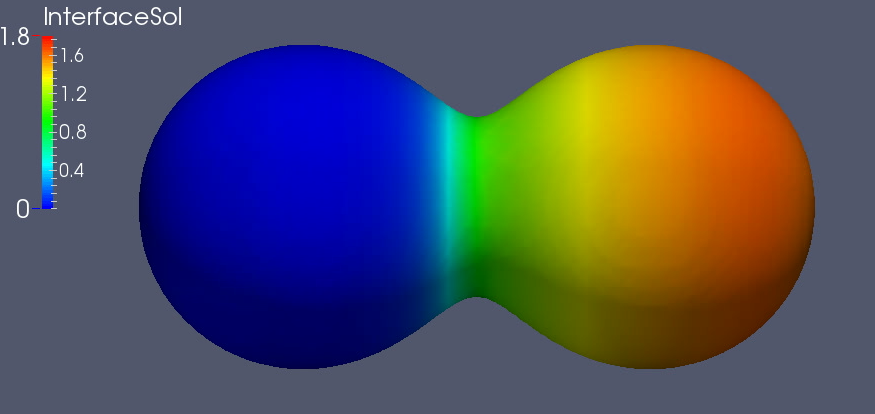}
 \end{minipage} \hfill
\begin{minipage}{0.48\textwidth}
  \includegraphics[width=\textwidth]{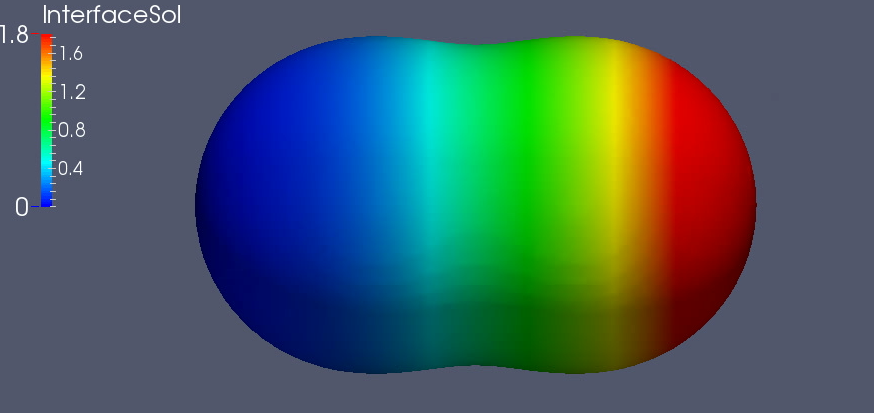}
  \end{minipage}\\[1ex]
\begin{minipage}{0.48\textwidth}
  \includegraphics[width=\textwidth]{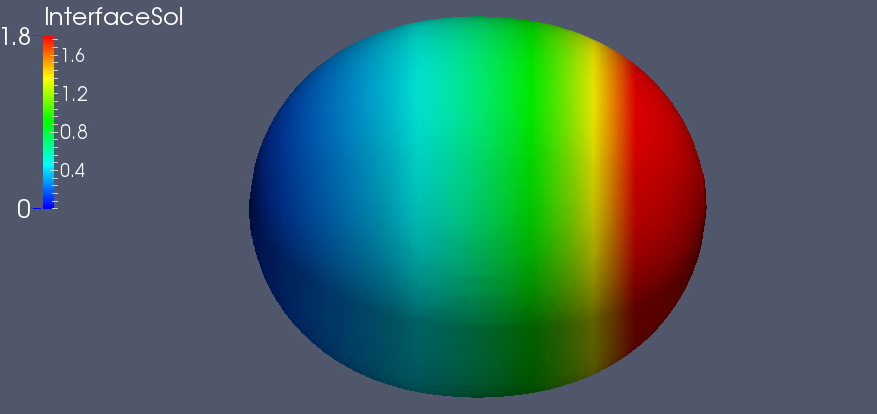}
 \end{minipage} \hfill
\begin{minipage}{0.48\textwidth}
  \includegraphics[width=\textwidth]{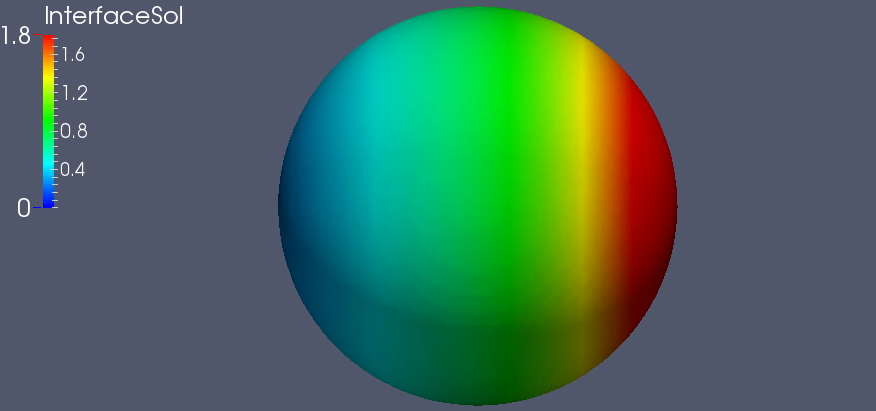}
  \end{minipage}
  \caption{\label{fig4} Snapshots of discrete solution, $l=4$, $\Delta t=\frac14$, at $t=0.25,0.5,0.75,1$.}
\end{centering}
\end{figure}

\begin{figure}[ht!]
  \centering
  \includegraphics[width=0.7\textwidth]{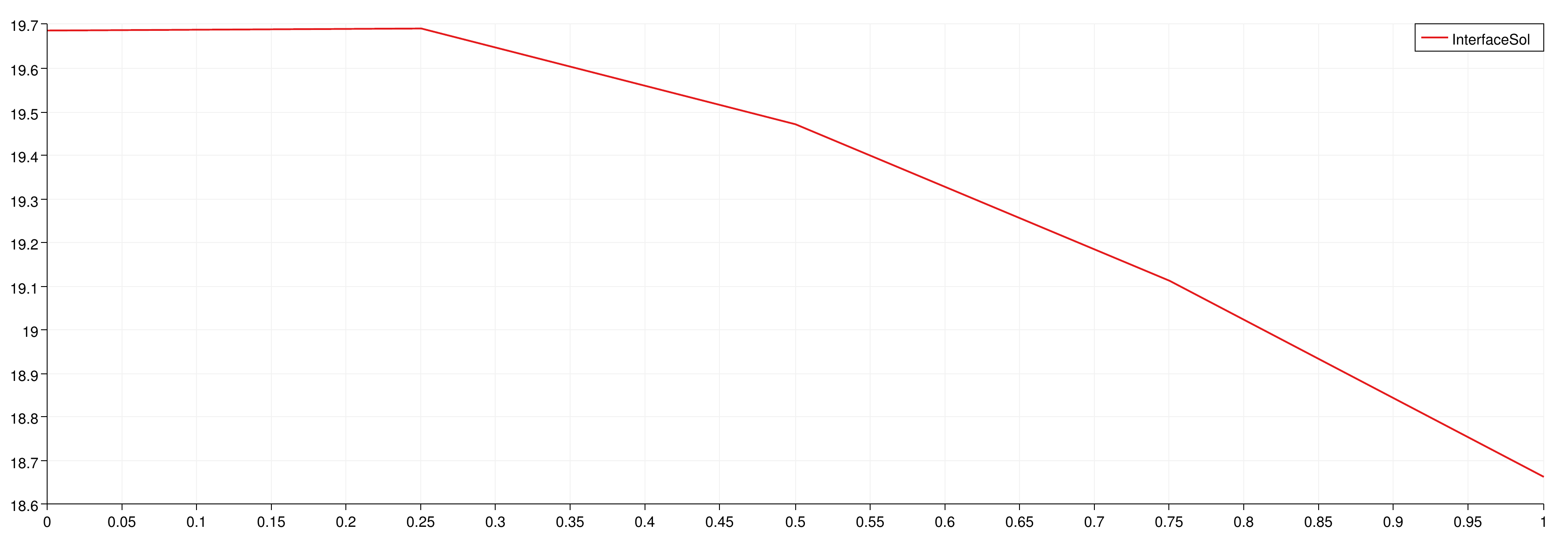}
  \caption{\label{fig:collision-mass-dt4}The total amount of surfactant $I_{4,1/4}$ over time.}
\end{figure}

 Finally, we show result for level $l=4$, but  with a large time step size $\Delta t=1/4$, i.e. we use use only four time steps to approximate the solution at $t=1$. The four discrete solutions at $t=0.25,0.5,0.75,1.0$ are shown in Figure~\ref{fig4}.  The total amount of surfactant is shown in Figure~\ref{fig:collision-mass-dt4}.


\section{DISCUSSION}

We presented a space-time finite element method for solving PDEs on evolving surfaces. The method is based on traces of outer  finite element spaces, is Eulerian in the sense that  $\Gamma(t)$ is not tracked by a mesh, and can easily be combined with both space and time adaptivity. No extension of the equation away from the surface is needed and thus the number of d.o.f. involved in
computations is optimal and comparable to methods in which $\Gamma$ is meshed directly. The computations are
done in a time-marching manner as common for parabolic equations.

The method has  second order  convergence in space and time  and conserves the mass in a weak sense, cf. \eqref{means}. In practice, an artificial mass flux can be experienced due to geometric errors resulting from  the approximation of $\Gamma(t)$.  In experiments,  the loss of mass was found to be  small and quickly vanishing if the mesh is refined.

The implicit definition of the surface evolution with the help of a level set function is  well suited for
numerical treatment of surfaces which undergo topological changes and experience singularities. This report
shows that the present space-time surface finite element method perfectly complements this property
and provides a robust technique for computing diffusion and transport along colliding surfaces.

\bibliographystyle{unsrt}
\bibliography{literatur}

\begin{thebibliography}{10}

\bibitem{DEreview}
G.~Dziuk and C.~M. Elliott.
\newblock Finite element methods for surface pdes.
\newblock {\em Acta Numerica}, 22:289--396, 2013.

\bibitem{Dziuk07}
G.~Dziuk and C.~Elliott.
\newblock Finite elements on evolving surfaces.
\newblock {\em IMA J. Numer. Anal.}, 27:262--292, 2007.

\bibitem{DziukElliot2013a}
G.~Dziuk and Ch.~M. Elliott.
\newblock {$L^2$-estimates} for the evolving surface finite element method.
\newblock {\em Mathematics of Computation}, 82:1--24, 2013.

\bibitem{AS03}
D.~Adalsteinsson and J.~A. Sethian.
\newblock Transport and diffusion of material quantities on propagating
  interfaces via level set methods.
\newblock {\em J. Comput. Phys.}, 185:271--288, 2003.

\bibitem{DziukElliot2010}
G.~Dziuk and C.~Elliott.
\newblock An eulerian approach to transport and diffusion on evolving implicit
  surfaces.
\newblock {\em Comput. Vis. Sci.}, 13:17–--28, 2010.

\bibitem{XuZh}
Jian-Jun Xu and Hong-Kai Zhao.
\newblock An {Eulerian} formulation for solving partial differential equations
  along a moving interface.
\newblock {\em Journal of Scientific Computing}, 19:573--594, 2003.

\bibitem{ORXsinum}
M.~Olshanskii, A.~Reusken, and X.~Xu.
\newblock An {Eulerian} space-time finite element method for diffusion problems
  on evolving surfaces.
\newblock NA\&SC Preprint~5, Department of Mathematics, University of Houston,
  2013.
\newblock Revised version submitted to SIAM J. Numer. Anal.

\bibitem{refJoerg}
J.~Grande.
\newblock Finite element methods for parabolic equations on moving surfaces.
\newblock IGPM Preprint 360, RWTH Aachen University, 2013.
\newblock Accepted for publication in SIAM J. Sci. Comput.

\bibitem{ORXsinum2}
M.~Olshanskii and A.~Reusken.
\newblock Error analysis of a space-time finite element method for solving
  {PDEs} on evolving surfaces.
\newblock IGPM Preprint 376, RWTH Aachen University, 2013.
\newblock Submitted.

\bibitem{OlshReusken08}
M.~Olshanskii, A.~Reusken, and J.~Grande.
\newblock A finite element method for elliptic equations on surfaces.
\newblock {\em SIAM J. Numer. Anal.}, 47:3339--3358, 2009.

\bibitem{OlshanskiiReusken08}
M.~Olshanskii and A.~Reusken.
\newblock A finite element method for surface {PDEs}: matrix properties.
\newblock {\em Numer. Math.}, 114:491--520, 2009.

\bibitem{DemlowOlshanskii12}
A.~Demlow and M.A. Olshanskii.
\newblock An adaptive surface finite element method based on volume meshes.
\newblock {\em SIAM J. Numer. Anal.}, 50:1624--1647, 2012.

\bibitem{James04}
A.J. James and J.~Lowengrub.
\newblock A surfactant-conserving volume-of-fluid method for interfacial flows
  with insoluble surfactant.
\newblock {\em J. Comp. Phys.}, 201(2):685--722, 2004.

\bibitem{GReusken2011}
S.~{Gro\ss} and A.~Reusken.
\newblock {\em Numerical Methods for Two-phase Incompressible Flows}.
\newblock Springer, Berlin, 2011.

\bibitem{TobiskaBook}
H.-G. Roos, M.~Stynes, and L.~Tobiska.
\newblock {\em Numerical Methods for Singularly Perturbed Differential
  Equations --- Convection-Diffusion and Flow Problems}, volume~24 of {\em
  Springer Series in Computational Mathematics}.
\newblock Springer-Verlag, Berlin, second edition, 2008.

\end{thebibliography}



\end{document}